\documentclass[11pt,reqno,a4paper]{amsart}

\usepackage[pdftex]{graphicx}
\usepackage{amsmath}
\usepackage{amssymb,amsthm,hyperref,caption}
\usepackage{color}
\definecolor{blau}{rgb}{0,0,0.75} 
\hypersetup{colorlinks,linkcolor=blau,citecolor=blue}

\allowdisplaybreaks

\newtheorem{theorem}{Theorem}

\theoremstyle{definition}

\newtheorem{example}{Example}

\newcommand{\N}{\ensuremath{\mathbb{N}}}

\newcommand{\R}{\ensuremath{\mathbb{R}}}

\newcommand{\ho}{\text{hook}}
\renewcommand{\deg}{\text{deg}}
\begin{document}

\author[M.~Kuba]{Markus Kuba}
\address{Markus Kuba\\
Institut f{\"u}r Diskrete Mathematik und Geometrie\\
Technische Universit\"at Wien\\
Wiedner Hauptstr. 8-10/104\\
1040 Wien -- HTL Wien 5 Spengergasse, Spengergasse 20, 1050 Wien, Austria} %
\email{kuba@dmg.tuwien.ac.at}

\title{A note on hook length formulas for trees}

\begin{abstract}
In this short note we discuss recent results on hook length formulas of trees unifying some earlier results, and explain hook length formulas naturally associated to families of increasingly labelled trees.
\end{abstract}


\maketitle

\section{Introduction}
In recent works Han~\cite{Han1,Han2} developed an expansion technique for deriving hook length formulas
for binary trees. Han's technique was applied by Chen et al.~\cite{Chen} to $k$-ary trees, 
plane trees, labelled trees and related families of forests. Moreover, Yang~\cite{Yang}
also gave extensions of Han's results for tree families related to $k$-ary trees and plane trees.
Furthermore, Sagan recently gave probabilistic proofs for some results of \cite{Han1,Yang,Chen}. 
In this short note we apply Han's expansion technique to families of weighted trees. 
The weights are similar to the weights in so-called simply generated tree families. This enables us to (re)-derive families of hook length formulas naturally associated to families of increasingly labelled trees, and to unify some results of~\cite{Chen,Yang}. First we give a quick review of the terminology. 
For a tree $T$ the hook length of a vertex $v\in T$, devoted by $h_v$, is the number of descendants 
of node $v$ with the convention that a node is counted as a descendant of itself. 
For a given weight function $\rho:\N^+\to\R$ we can associate to any given tree $T$ of some tree family, e.g., $T$ being a binary tree or a plane tree, a weight 
\begin{equation}
\label{Hookeqn1}
w_{\ho}(T)=\prod_{v\in T}\rho(h_v).
\end{equation}
Han~\cite{Han1} derived the following result for
binary trees, where $\mathcal{B}(n)$ denotes the set of binary trees of size $n$, i.e.~with $n$ vertices.
If the relation 
\begin{equation}
\label{HookHan1}
\sum_{n\ge 1}\Big(\sum_{T\in \mathcal{B}(n)}w_{\ho}(T)\Big)z^n=F(z),
\end{equation}
holds, then the weight function $\rho$ satisfies
\begin{equation}
\label{HookHan2}
\rho(n)=\frac{[z^n]F(z)}{[z^{n-1}](1+F(z))^2}.
\end{equation}
Here $[z^n]$ denotes the extraction of coefficient operator, $[z^n]F(z)$ denotes the coefficient
of $z^n$ in the formal power series expansion of $F(z)$.
The result~\ref{HookHan1},~\ref{HookHan2} was extended by Chen et al.~\cite{Chen} and Yang~\cite{Yang}
using Han's method, and many hook length formulas were deduced. We extend the relation~\ref{HookHan1},~\ref{HookHan2} to trees which are weighted in analogy to simply generated tree families $\mathcal{T}$, sometimes also called simple varieties of trees, Flajolet and Sedgewick~\cite{FlaSed2009}. Important steps in this direction where already carried out by Chen et al.~\cite{Chen} and Yang~\cite{Yang}.

\smallskip

A class $\mathcal{T}$ of simply generated trees can be defined in
the following way. A sequence of non-negative real numbers
$(\varphi_{k})_{k \ge 0}$ with $\varphi_{0} > 0$ ($\varphi_{k}$
can be seen as the multiplicative weight of a node with out-degree
$k$) is used to define the weight $w(T)$ of any \emph{ordered tree} $T$
by $w_{\deg}(T) := \prod_{v} \varphi_{d(v)}$, where $v$ ranges over all
vertices of $T$ and $d(v)$ is the out-degree (the number of
children) of $v$ (in order to avoid degenerate cases we always
assume that there exists a $k \ge 2$ such that $\varphi_{k} > 0$).
The family $\mathcal{T}$ consists then of all trees $T$ with $w_{\deg}(T)
\neq 0$ together with their weights $w_{\deg}(T)$. It follows further
that for a given degree-weight sequence $(\varphi_{k})_{k \ge 0}$
the generating function $T(z) := \sum_{n \ge 1} T_{n} z^{n}$ of
the quantity total weights $T_{n} := \sum_{|T|=n} w_{\deg}(T)$, where
$|T|$ denotes the size of the tree $T$, satisfies the functional
equation
\begin{equation}
   \label{eqnb1}
   T(z) = z \varphi\big(T(z)\big),
\end{equation}
where the degree-weight generating function $\varphi(t)$ is given
by $\varphi(t) = \sum_{k \ge 0} \varphi_{k} t^{k}$. 
By suitable choices for the degree generating function $\varphi(t)$ one obtains several important tree families as special cases,
such as binary trees $\varphi(t)=(1+t)^2$, $k$-ary trees $\varphi(t)=(1+t)^k$, plane trees $\varphi(t)=1/(1-t)$, labelled trees, 
$\varphi(t)=e^t$, see~\cite{FlaSed2009}.

\section{Main result}
\begin{theorem}
\label{HookThemain}
For a given degree weight generating function $\varphi(t)$ assume that the relation 
\begin{equation}
\sum_{n\ge 1}\Big(\sum_{T\in \mathcal{T}(n)}w_{\deg}(T)w_{\ho}(T)\Big)z^n=F(z),
\end{equation}
holds, where $w_{\deg}(T) := \prod_{v} \varphi_{d(v)}$, $w_{\ho}(T)=\prod_{v\in T}\rho(h_v)$, and $\mathcal{T}$ is the family of ordered trees. Then the weight function $\rho$ satisfies
\begin{equation}
\rho(n)=\frac{[z^n]F(z)}{[z^{n-1}]\varphi(F(z))},\quad n\ge 1.
\end{equation}
\end{theorem}

\begin{example}
One obtains Theorem~2.1, 3.1, 4.1 of Chen et al.~\cite{Chen} by considering $\varphi(t)=(1+t)^k$, $\varphi(t)=1/(1-t)$ and $\varphi(t)=e^{t}$.
Moreover, the families of degree weight generating functions considered by Yang~\cite{Yang} are
given by $\varphi(t)=(1+s\cdot t)^m$, for real parameters $s$ and $m$.
\end{example}

\begin{proof}
Let $F_n=[z^n]F(z)$. By definition we have
\begin{equation*}
F_n=\sum_{T\in \mathcal{T}(n)}w_{\deg}(T)w_{\ho}(T).
\end{equation*}
Following Han~\cite{Han1}, see also Chen et al.~\cite{Chen} and Yang~\cite{Yang}, we use the so-called top-bottom decomposition of a tree, decomposing a tree into a root $r$ 
and a subtrees attached to the root. Assuming that the root of a given size $n$ tree $T$ has outdegree $j\ge 1$. Then we get a degree weight factor $\varphi_{d(r)}= \varphi_j$, a hook length weight factor $\rho(h_{r}) =\rho(n)$, and factors corresponding to $j$ non-empty subtrees $T_{\ell}$, $1\le \ell\le j$, dangling from the root,
\begin{equation*}
w_{\deg}(T)w_{\ho}(T)=\varphi_{j}\cdot \rho(n) \cdot \prod_{\ell=1}^{j}\Big(w_{\deg}(T_{\ell})w_{\ho}(T_{\ell})\Big).
\end{equation*}
By considering all 
possible ways to build a tree of size $n$ we get
\begin{equation*}
\sum_{T\in \mathcal{T}(n)}w_{\deg}(T)w_{\ho}(T)
=  \rho(n)\sum_{j\ge 1}\varphi_j \sum_{\sum_{\ell=1}^{n}n_{\ell}=n-1}\prod_{i=1}^{j}F_{n_{\ell}}
= \rho(n)[z^{n-1}]\varphi\big(F(z)\big).
\end{equation*}
One still has to check the initial case $n=1$. In accordance to the definition $F_1=\varphi_0\rho(1)$ we get for $n=1$ the equation $\rho(1)=\frac{[z^1]F(z)}{[z^{0}]\varphi(F(z))}=\frac{F_1}{\varphi_0}$, which proves the stated result.
\end{proof}

One may also derive expansions of $\varphi(F(z))$, related to forests, as shown below.
\begin{theorem}
For a given degree weight generating function $\varphi(t)$ assume that the relation 
\begin{equation}
\begin{split}
&\varphi_0+\sum_{n\ge 1}\Big(\sum_{T\in \mathcal{F}(n)}w_{\deg}(T)w_{\ho}(T)\Big)z^n = G(z)\\
&=\varphi\bigg(\sum_{n\ge 1}\Big(\sum_{T\in \mathcal{T}(n)}w_{\deg}(T)w_{\ho}(T)\Big)z^n\bigg)=\varphi(F(z)),
\end{split}
\end{equation}
holds, for a class of forests $\mathcal{F}$, where $w_{\deg}(T) := \prod_{v} \varphi_{d(v)}$, $w_{\ho}(T)=\prod_{v\in T}\rho(h_v)$, and $\mathcal{T}$ is the family of ordered trees. Then the weight function $\rho$ satisfies 
\begin{equation}
\rho(n)=\frac{[z^n]\varphi^{[-1]}\big(G(z)\big)}{[z^{n-1}]G(z)},\quad n\ge 1,
\end{equation}
where $\varphi^{[-1]}(t)$ denotes the inverse degree weight generating function, $\varphi(\varphi^{[-1]}(t))=t$.
\end{theorem}
\begin{proof}
By application of the inverse function $\varphi^{[-1]}(t)$ we obtain the equation
\begin{equation*}
F(z)=\sum_{n\ge 1}\Big( \sum_{T\in \mathcal{T}(n)}w_{\deg}(T)w_{\ho}(T)\Big)z^n=\varphi^{[-1]}\big(G(z)\big).
\end{equation*}
Hence, by Theorem~\ref{HookThemain} we obtain 
\begin{equation}
\rho(n)=\frac{[z^n]\varphi^{[-1]}\Big(G(z)\Big)}{[z^{n-1}]\varphi\Big(\varphi^{[-1]}\big(G(z)\big)\Big)}
=\frac{[z^n]\varphi^{[-1]}\Big(G(z)\Big)}{[z^{n-1}]G(z)},\quad n\ge 1.
\end{equation}
\end{proof}
\begin{example}
One readily obtains Theorem~3.4 and 4.4 of Chen et al.~\cite{Chen} concerning plane forests, and forests of labelled trees by considering 
the degree weight generating functions $\varphi(t)=1/(1-t)$, 
such that 
\begin{equation*}
\varphi^{[-1]}(t)=1-\frac{1}{t},\quad  \rho(n)=\frac{[z^n]\varphi^{[-1]}\big(G(z)\big)}{[z^{n-1}]G(z)}= \frac{-[z^n]\big(G(t)\big)^{-1}}{[z^{n-1}]G(z)},
\end{equation*}
and $\varphi(t)=e^{t}$, such that
\begin{equation*}
\varphi^{[-1]}(t)=\log(t),\quad  \rho(n)=\frac{[z^n]\varphi^{[-1]}\big(G(z)\big)}{[z^{n-1}]G(z)}= \frac{[z^n]\log G(t)}{[z^{n-1}]G(z)}.
\end{equation*}
\end{example}

\section{Increasing trees and hook length formulas}
\subsection{Definition}
Increasing trees are labelled trees where
the nodes of a tree of size $n$ are labelled by distinct integers of the set
$\{1, \dots, n\}$ in such a way that each sequence of labels along any branch starting
at the root is increasing. Formally, a class $\mathcal{T}$ of a simple family of increasing trees can be defined 
in the following way, similar to simply generated trees. A sequence of non-negative numbers $(\varphi_{k})_{k \ge 0}$
with $\varphi_{0} > 0$ is used to define the weight $w(T)$ of any ordered tree $T$ by
$w_{\deg}(T)=\prod_{v} \varphi_{d(v)}$, where $v$ ranges over all vertices of $T$ and $d(v)$ is the
out-degree of $v$ (we always assume that it exists a $k \ge 2$ with $\varphi_{k} > 0$).
Furthermore, $\mathcal{L}(T)$ denotes the set of different increasing
labellings of the tree $T$ with distinct integers $\{1, 2, \dots, |T|\}$,
where $|T|$ denotes the size of the tree $T$, and $L(T) :=
\big|\mathcal{L}(T)\big|$ its cardinality.
Then the family $\mathcal{T}$ consists of all trees $T$ together with their
weights $w_{\deg}(T)$ and the set of increasing labellings $\mathcal{L}(T)$.
For a given degree-weight sequence $(\varphi_{k})_{k \ge 0}$
with a degree-weight generating function
$\varphi(t) := \sum_{k \ge 0} \varphi_{k} t^{k}$, we define now the total weights
by $T_{n} := \sum_{|T|=n} w_{\deg}(T) \cdot L(T)$.
It follows then that the exponential generating function
$T(z) := \sum_{n \ge 1} T_{n} \frac{z^{n}}{n!}$
satisfies the autonomous first order differential equation
\begin{equation}
   \label{eqnz1}
   T'(z) = \varphi\big(T(z)\big), \quad T(0)=0.
\end{equation}

\subsection{Natural hook length formulas}
For any unlabelled rooted tree $T$ having $n$ vertices the following formula for $L(T) :=
\big|\mathcal{L}(T)\big|$ is well known.
\begin{equation*}
L(T)=\frac{n!}{\prod_{v\in T}h_{v}}.
\end{equation*}
Hence, the definition of the total weights $T_n$ naturally involves hook lengths,
\begin{equation*}
\sum_{T\in \mathcal{T}(n)} w_{\deg}(T) \cdot \frac{n!}{\prod_{v\in T}h_{v}} = T_n.
\end{equation*}
\begin{example}
Let $\varphi(t)=1/(1-t)^{\alpha}$, $\alpha>0$. Then, according to~\eqref{eqnz1}
$T'(z)=(1-T(z))^{-\alpha}$, such that $T(z)=1-(1-(\alpha+1)z)^{1/(\alpha+1)}$, 
and $T_n=[z^n]T(z)=(\alpha+1)^{n-1}(n-1)!\binom{n-1-\frac{1}{\alpha+1}}{n-1}$. 
Hence, we obtain 
\begin{equation*}
\sum_{T\in \mathcal{T}(n)} 
\Big(\prod_{v} \binom{\alpha -1+ d(v)}{d(v)}\Big) \cdot \frac{n!}{\prod_{v\in T}h_{v}} = (\alpha+1)^{n-1}(n-1)!\binom{n-1-\frac{1}{\alpha+1}}{n-1}.
\end{equation*}
In order to derive the formula above from Theorem~\ref{HookThemain} 
one adapts $\varphi(t)=1/(1-t)^{\alpha}$ and $F(z)=1-(1-(\alpha+1)z)^{1/(\alpha+1)}$ 
such that 
\begin{equation*}
\rho(n)=\frac{[z^n]F(z)}{[z^{n-1}]\varphi(F(z))}=\frac{[z^n]F(z)}{[z^{n-1}]F'(z)} = \frac{1}{n},
\end{equation*}
which also leads to the stated formula.  
\end{example}
The case $\alpha=1$ appears as Example 3.2 in~\cite{Chen} and was the primary motivation/inspiration for writing this short note.
Similarly, related formulas for $k$-ary increasing trees, etc., can be obtained from Theorem~\ref{HookThemain} .

\end{document}